\documentclass[reqno,12pt]{amsart}
\usepackage{amsfonts}
\usepackage{amssymb}
\usepackage{bbm}
\usepackage{times}
\usepackage{bbm}
\usepackage{amsmath}
\usepackage{txfonts}
\usepackage{stmaryrd}
\usepackage{mathrsfs}
\usepackage{setspace}

\title{On Invariant subspace for hyponormal operators}

\author{Junfeng Liu}

\address{Junfeng Liu\\ Department of Mathematics, School of Science, Zhejiang University of Science and Technology, Hangzhou, 310023, P. R. China. }
\email{jfliu997@sina.com}

\date{}
\newtheorem{theorem}{Theorem}

\newtheorem{corollary}{Corollary}

\newtheorem{case}{Case}

\begin{document}

\maketitle

\begin{abstract}
There is a resent paper claiming that every hyponormal operator which is not a multiple of
the identity (operator) has a nontrivial hyperinvariant subspace. If this claim is true, then every hyponormal operator has a nontrivial
invariant subspace, and every subnormal operator which is not multiple of the identity has a nontrivial
hyperinvariant subspace. But we find out that the proof of the above claim go wrong.
Therefore the invariant subspace problem of the hyponormal operator and the hyperinvariant problem of subnormal
operator remains unresolved. Moreover, it is well known that these are two important research topics in operator
theory that have not been solved for a long time, and that a lot of people have been trying to solved the two problems.
So it is very meaningful to clarify whether these two problems have been solved.
\end{abstract}

\textbf{Keywords:}  invariant subspace; hyperinvariant subspace;  subnormal operator; hyponormal operator; Hilbert space.

\textbf{2020 MSC:} 47B20; 47A15.

\section{Introduction and Preliminary}

It is well known that the invariant subspace problem is one of the famous unsolved problem in operator theory
(cf. \cite{na54}, \cite{ab49}, \cite{bro87}, \cite{hpr82}, \cite{jl192}, \cite{sm22}, \cite{cjr85} and their references),
and that there are many famous mathematical journals that published research papers
on invariant subspaces (see\cite{na54}, \cite{ab49}, \cite{bro87}, \cite{pe87} and so on). This problem is as follows:
does every bounded linear operator on an infinite dimensional Hilbert space have a nontrivial invariant subspace
(see \cite{hpr82} p.100, \cite{jl20}, \cite{jl22}, \cite{sm22} and so on).
It is also well known that if an operator has a nontrivial hyperinvariant subspace then it has a nontrivial
invariant subspace, and that every subnormal operator is a hyponormal operator, but the inverse
of there two propositions does not hold.

In 1978, S. Brown \cite{sb78} proved every subnormal operator has a nontrivial invariant
subspace. But the hyperinvariant subspace problem of subnormal operator has not been solved for a long time
(see \cite{bea88} p.155 and so on).

In 1987, S. Brown published a paper in Ann. of Math., in which he proved that every hyponormal operator
with thick spectrum has a nontrivial invariant subspace.
But the invariant subspace problem of the Hyponormal operator is still open (cf. \cite{sm22} p.234 and so on). Since 1987, a lot
of people (including ourself) have been trying to prove every (general) hyponormal
operator has a nontrivial invariant subspace. These facts show that the
invariant subspace problem of hyponormal operator and the hyperinvariant subspace of the subnormal
operator are very important and very meaningful
research topics.

\section{Main Results}

In 2022, S. Mecheri published a paper \cite{sm22} in which he given the following theorem:

\begin{theorem}{Let $H$ be a Hilbert space and let $T$ be a bounded linear operator on $H$.
Assume that $T$ is not a multiple of the identity operator $I$ and unitary equivalent
to an upper triangular operator matrix on $H$, say
\begin{eqnarray}\label{23t1}
T=\left(
\begin{array}{cccccccccc}
T_{11}&T_{12}&T_{13}&T_{14}&T_{15}&\cdots&T_{1n}\\
0&T_{22}&T_{23}&T_{24}&T_{25}&\cdots&T_{2n}\\
0&0&T_{33}&T_{34}&T_{35}&\cdots&T_{3n}\\
0&0&0&T_{44}&T_{45}&\cdots&T_{4n}\\
0&0&0&0&T_{55}&\cdots&T_{5n}\\
\vdots&\vdots&\vdots&\ddots&\vdots&\vdots&\vdots&\\
0&0&0&0&0&\cdots&T_{nn}\\
\end{array}
\right).
\end{eqnarray}
If $T_{11}$ and $T_{nn}$ are M-hyponormal operators, then $T$ has a nontrivial hyperinvariant subspace.
}\end{theorem}

This is just Theorem 3.8 in \cite{sm22}. It is well known that every subnormal operator is
a hyponormal operator and that every hyponormal operator is a M-hyponormal operator.
Thus by the above theorem, the another given the following corollary:

\begin{corollary}{
Every M-hyponormal operator which is not a multiple of the identity has a nontrivial hyperinvariant
subspace. In particular, every hyponormal operator
(containing subnormal operator) that is not a multiple of the identity has a nontrivial hyperinvariant
subspace.
}\end{corollary}

This is just Corollary 3.12 in \cite{sm22}. The author of \cite{sm22} divided the proof of the above
theorem into two case as follows:

\begin{case}\label{23c1}
{The spectrum of $T_{11}$ and $T_{nn}$ consists of more than one point.
}\end{case}

\begin{case}\label{23c2}{
The spectrum of $T_{11}$ and $T_{nn}$ consist of only one point.
}\end{case}
But we find out that there are an error and a query in the proof of Case \ref{23c1} and Case \ref{23c2}
respectively. The details as follows:

\medskip

\textbf{Error 1.} In Case \ref{23c1}, the author
of \cite{sm22} written that "if the spectrum of $T_{11}$ and $T_{nn}$ consists of more than one point.
$\cdots$, where $N$ is normal and $TN=NT$. Let $N=\int_{\sigma(N)}zdE(z)$ be spectral decomposition of
the normal operator $N$. $\cdots$. Assume that $\{B_1,B_2\}$ is a partition of $\sigma(T_{11})$ and $\{B_3,B_4\}$
is a partition of $T_{nn}$ into Borel sets, where $B_1$ and $B_4$ are disjoint closed sets.
Let $T_{11}^0=T_{11}|_{E(B_1)H}, ~T_{11}^1=T_{11}|_{E(B_2)H},$
$T_{nn}^0=T_{nn}|_{E(B_3)H}$ and $T_{nn}^1=T_{nn}|_{E(B_4)H}$. Then $T_{11}=T_{11}^0\oplus T_{11}^1$ and
$T_{nn}=T_{nn}^0\oplus T_{nn}^1$."
According to the spectral decomposition theorem of the normal operator $N$,
where $E(B)$ is the spectral measure of the Borel subset $B$ of $\sigma(N)$, and is a projection operator
from $H$ onto $E(B)H$. (cf. \cite{jbc90} p.256-265).

Now, we give a counterexample to illustrate that the above claims $T_{11}=T_{11}^0\oplus T_{11}^1$
and $T_{nn}=T_{nn}^0\oplus T_{nn}^1$ are not true.

\medskip

\textbf{Counterexample 1.} In Theorem 3.8 of \cite{sm22}, the operator $T$ is an upper triangular
operator matrix on the Hilbert space $H$. Now, let $H=l_2\oplus l_2\oplus \cdots \oplus l_2$
is the orthogonal sum of $n$ copies of $l_2$. Let $N$ be an operator on $H$ whose representation matrix is as
follows:
\begin{eqnarray}\label{23t2}
\left(
\begin{array}{cccc}
1&0&0&\cdots\\
0&0&0&\cdots\\
0&0&0&\cdots\\
\vdots&\vdots&\vdots&\ddots
\end{array}
\right).
\end{eqnarray}
Let $T_{11}$ and $T_{nn}$ be the operators on $l_2$ with the representation matrix

\begin{eqnarray}\label{23t3}
\left(
\begin{array}{cccc}
0&0&0&\cdots\\
0&2&0&\cdots\\
0&0&0&\cdots\\
\vdots&\vdots&\vdots&\ddots
\end{array}
\right).
\end{eqnarray}

(a). First of all, we show that $N$ is a normal operator on $H$. In fact, it is clear that for each
$x=(x_1,x_2,\cdots)\in H$, we have
\begin{eqnarray}\label{23t4}
Nx=\left(
\begin{array}{cccc}
1&0&0&\cdots\\
0&0&0&\cdots\\
0&0&0&\cdots\\
\vdots&\vdots&\vdots&\ddots
\end{array}
\right)
\left(\begin{array}{c}
x_1\\
x_2\\
x_3\\
\vdots
\end{array}
\right)
=\left(\begin{array}{c}
x_1\\
0\\
0\\
\vdots
\end{array}
\right).
\end{eqnarray}
It is easy to see that $N$ is a bounded linear operator, and that the representation matrix of the adjoint
operator $N^*$ of $N$ is also the matrix in (\ref{23t2}). Indeed, for each $y=(y_1,y_2,\cdots)\in H$,
we have
\begin{eqnarray*}
N^*y=\left(
\begin{array}{cccc}
1&0&0&\cdots\\
0&0&0&\cdots\\
0&0&0&\cdots\\
\vdots&\vdots&\vdots&\ddots
\end{array}
\right)
\left(\begin{array}{c}
y_1\\
y_2\\
y_3\\
\vdots
\end{array}
\right)
=\left(\begin{array}{c}
y_1\\
0\\
0\\
\vdots
\end{array}
\right).
\end{eqnarray*}
Consequently, $\langle Nx,y\rangle=x_1y_1=\langle x,N^*y\rangle$. Therefore, $N^*$ is the adjoint operator
of $N$. Moreover, we have
\begin{eqnarray*}
NN^*x=\left(
\begin{array}{cccc}
1&0&0&\cdots\\
0&0&0&\cdots\\
0&0&0&\cdots\\
\vdots&\vdots&\vdots&\ddots
\end{array}
\right)
\left(\begin{array}{c}
x_1\\
0\\
0\\
\vdots
\end{array}
\right)
=\left(\begin{array}{c}
x_1\\
0\\
0\\
\vdots
\end{array}
\right).
\end{eqnarray*}
and
\begin{eqnarray*}
N^*Nx=\left(
\begin{array}{cccc}
1&0&0&\cdots\\
0&0&0&\cdots\\
0&0&0&\cdots\\
\vdots&\vdots&\vdots&\ddots
\end{array}
\right)
\left(\begin{array}{c}
x_1\\
0\\
0\\
\vdots
\end{array}
\right)
=\left(\begin{array}{c}
x_1\\
0\\
0\\
\vdots
\end{array}
\right).
\end{eqnarray*}
Therefore $NN^*=N^*N$. Hence $N$ is normal.

(b). Next, we show that $\sigma(N)=\{0,1\}$. In fact, when $\lambda=0$, by (\ref{23t4}) we have
$$(N-\lambda I)x=Nx=(x_1,0,0\cdots),$$
where $I$ denotes the identity.
So $N-\lambda I$ is not surjective from $H$ onto itself. This implies that $\lambda=0\in \sigma(N)$.
When $\lambda=1$, we have
\begin{eqnarray*}
(N-\lambda I)x=\left(
\begin{array}{cccc}
0&0&0&\cdots\\
0&-1&0&\cdots\\
0&0&-1&\cdots\\
\vdots&\vdots&\vdots&\ddots
\end{array}
\right)
\left(\begin{array}{c}
x_1\\
x_2\\
x_3\\
\vdots
\end{array}
\right)
=\left(\begin{array}{c}
0\\
-x_2\\
-x_3\\
\vdots
\end{array}
\right).
\end{eqnarray*}
Again, $N-\lambda I$ is not surjective from $H$ onto itself, so that $\lambda=1\in \sigma(N)$.

When $\lambda\neq 0$ and $\lambda \neq 1$ , for $x=(x_1, x_2, \cdots)\in H$, put
\begin{eqnarray}\label{23t5}
(N-\lambda I)x=\left(
\begin{array}{cccc}
1-\lambda&0&0&\cdots\\
0&-\lambda&0&\cdots\\
0&0&-\lambda&\cdots\\
\vdots&\vdots&\vdots&\ddots
\end{array}
\right)
\left(\begin{array}{c}
x_1\\
x_2\\
x_3\\
\vdots
\end{array}
\right)
=\left(\begin{array}{c}
(1-\lambda)x_1\\
-\lambda x_2\\
-\lambda x_3\\
\vdots
\end{array}
\right)
=
\left(\begin{array}{c}
y_1\\
y_2\\
y_3\\
\vdots
\end{array}
\right)=y,
\end{eqnarray}
where $y=(y_1,y_2,\cdots)\in H$.
In (\ref{23t5}), set $y=(y_1,y_2,\cdots)=(0,0,\cdots)$, then
$x=(x_1,x_2,\cdots)=(0,0,\cdots)$. Therefore $N-\lambda I$ is injective from $H$ into itself.
Moreover, it is easy to see from (\ref{23t5}) that $N-\lambda I$ is surjective from $H$
onto itself whenever $\lambda\neq0$ and $\lambda\neq 1$. Consequently, $N-\lambda I$ is a
bijective from $H$ onto itself whenever $\lambda\neq0$ and $\lambda\neq 1$.
Hence $\lambda\in \rho(N)=\mathbb{C} \setminus\sigma(N)$
whenever $\lambda\neq0$ and $\lambda\neq 1$. Thus we conclude that $\sigma(N)=\{0,1\}$.

(c). we now show that $T_{11}$ is an M-hyponormal operator on $l_2$. In fact, it is clear that for each
$x=(x_1,x_2,\cdots)\in l_2$, we have
\begin{eqnarray}\label{23t6}
T_{11}x=\left(
\begin{array}{cccc}
0&0&0&\cdots\\
0&2&0&\cdots\\
0&0&0&\cdots\\
\vdots&\vdots&\vdots&\ddots
\end{array}
\right)
\left(\begin{array}{c}
x_1\\
x_2\\
x_3\\
\vdots
\end{array}
\right)
=\left(\begin{array}{c}
0\\
2x_2\\
0\\
\vdots
\end{array}
\right).
\end{eqnarray}
As in Part (a), it is easy to see that the representation matrix of the adjoint operator
$T^*_{11}$ of $T_{11}$ is also the matrix in (\ref{23t3}). Thus, for each $y=(y_1,y_2,\cdots)\in l_2$,
we have
\begin{eqnarray*}
T^*_{11}y=\left(
\begin{array}{cccc}
0&0&0&\cdots\\
0&2&0&\cdots\\
0&0&0&\cdots\\
\vdots&\vdots&\vdots&\ddots
\end{array}
\right)
\left(\begin{array}{c}
y_1\\
y_2\\
y_3\\
\vdots
\end{array}
\right)
=\left(\begin{array}{c}
0\\
2y_2\\
0\\
\vdots
\end{array}
\right).
\end{eqnarray*}
As in Part (a), it can be shown that
$T_{11}T^*_{11}=T^*_{11}T_{11}$. Therefore $T_{11}$ is normal. Hence $T_{11}$ is M-hyponormal.
Similarly, $T_{nn}$ is  M-hyponormal.

(d). Moreover, as in Part (b), it can be shown that $\sigma(T_{11})=\{0,2\}$, $\sigma(T_{nn})=\{0,2\}$

(e). It is easy to see from the definition of those operators that for any
$x=(x_{11}, x_{12}, \cdots,x_{n1}, x_{n2},\cdots)\in H$, we have
\begin{eqnarray}\label{23t7}
Tx=\left(
\begin{array}{cccccccc}
0&0&0&\cdots&0&0&0&\cdots\\
0&2&0&\cdots&0&0&0&\cdots\\
0&0&0&\cdots&0&0&0&\cdots\\
\vdots&\vdots&\vdots&\ddots&\vdots&\vdots&\vdots&\ddots\\
0&0&0&\cdots&0&0&0&\cdots\\
0&0&0&\cdots&0&2&0&\cdots\\
0&0&0&\cdots&0&0&0&\cdots\\
\vdots&\vdots&\vdots&\ddots&\vdots&\vdots&\vdots&\ddots\\
\end{array}
\right)
\left(\begin{array}{c}
x_{11}\\
x_{12}\\
x_{13}\\
\vdots\\
x_{n1}\\
x_{n2}\\
x_{n3}\\
\vdots
\end{array}
\right)
=\left(\begin{array}{c}
0\\
2x_{12}\\
0\\
\vdots\\
0\\
2x_{n2}\\
0\\
\vdots
\end{array}
\right).
\end{eqnarray}
Now, it is easy to see that from (\ref{23t4}) and (\ref{23t7}) that
\begin{eqnarray*}
TNx=\left(
\begin{array}{cccccccc}
0&0&0&\cdots&0&0&0&\cdots\\
0&2&0&\cdots&0&0&0&\cdots\\
0&0&0&\cdots&0&0&0&\cdots\\
\vdots&\vdots&\vdots&\ddots&\vdots&\vdots&\vdots&\ddots\\
0&0&0&\cdots&0&0&0&\cdots\\
0&0&0&\cdots&0&2&0&\cdots\\
0&0&0&\cdots&0&0&0&\cdots\\
\vdots&\vdots&\vdots&\ddots&\vdots&\vdots&\vdots&\ddots\\
\end{array}
\right)
\left(\begin{array}{c}
x_1\\
0\\
0\\
\vdots\\
0\\
0\\
0\\
\vdots
\end{array}
\right)
=\left(\begin{array}{c}
0\\
0\\
0\\
\vdots\\
0\\
0\\
0\\
\vdots
\end{array}
\right).
\end{eqnarray*}
and
\begin{eqnarray*}
NTx=\left(
\begin{array}{cccccccc}
1&0&0&\cdots&0&0&0&\cdots\\
0&0&0&\cdots&0&0&0&\cdots\\
0&0&0&\cdots&0&0&0&\cdots\\
\vdots&\vdots&\vdots&\ddots&\vdots&\vdots&\vdots&\ddots\\
0&0&0&\cdots&0&0&0&\cdots\\
0&0&0&\cdots&0&0&0&\cdots\\
0&0&0&\cdots&0&0&0&\cdots\\
\vdots&\vdots&\vdots&\ddots&\vdots&\vdots&\vdots&\ddots\\
\end{array}
\right)
\left(\begin{array}{c}
0\\
2x_{12}\\
0\\
\vdots\\
0\\
2x_{n2}\\
0\\
\vdots
\end{array}
\right)
=\left(\begin{array}{c}
0\\
0\\
0\\
\vdots\\
0\\
0\\
0\\
\vdots
\end{array}
\right).
\end{eqnarray*}
Thus $TN=NT$.

Now we can conclude that $N$ is a normal operator on $H$ and $T_{11}$, $T_{nn}$
are M-hyponormal operator with $TN=NT$, and that $\sigma(N)=\{0,1\}$ and $\sigma(T_{11})=\sigma(T_{nn})=\{0,2\}$.
Consequently, $T,~N,~T_{11}$ and $T_{nn}$ satisfy all conditions in the proof of Case \ref{23c1} of Theorem 3.8 of \cite{sm22}.
Now, let $\{B_1,B_2\}$ be a partition of $\sigma(T_{11})$ into Borel subsets, then $B_1=\{0\}$, $B_2=\{2\}$ (or
$B_1=\{2\}$, $B_2=\{0\}$). It is clear that $B_2=\{2\}$ is not the (Borel) subset of $\sigma(N)$.
Thus by the spectral theorem of the normal operator (\cite{jbc90} p.263), the operator $E(B_2)$ is not
defined, so that the operator $T^1_{11}=T_{11}|_{E(B_2)H}$ is not defined. Thus the claim $T_{11}=T^0_{11}\oplus T^1_{11}$
is not true. Similarly, the claim $T_{nn}=T_{nn}^0\oplus T^1_{nn}$ is also not true.

\medskip

\textbf{Note 1.} In the above discussion, as usual, the spectral measure of the normal operator $N$
is defined on the Borel subsets of $\sigma(N)$. Moreover, in \cite{br97} on p.6, the spectral measure
$E$ of the spectral operator may be defined on the Borel subsets of $\mathbb{C}$ with $E(\phi)=0$
and $E(\mathbb{C})=I$, where $\phi$ denotes the empty set, and $I$ denotes the identity.
But it follows from the properties of the spectral integral of the normal operator $N$ that when $\phi_{_0}(z)\equiv 1$,
$$I=\phi_{_0}(N)=\int_{\sigma(N)}\phi_{_0}(z)dE(z)=\int_{\sigma(N)}dE(z)$$
(cf. \cite{jbc90} p.264 Theorem 2.3 and so on). Next, by the definition of the spectral integral, we have
$$\int_{\sigma(N)}dE(z)=E(\sigma(N)).$$
(cf. \cite{jbc90} p.258 proposition 1.10 and so on). Consequently
$$E(\sigma(N))=I.$$
In this case, for any Borel subset $B$ of $\mathbb{C}$, if $B\cap \sigma(N)=\phi$, the $E(B)=0$.
Consequently, in the above discussion, $\sigma(N)=\{0,1\}$, $B_2=\{2\}$, $\sigma(B_2\cap \sigma(N))=\phi$,
$E(B_2)=0$. Therefore $E(B_2)H=\{0\}$. Hence $T^1_{11}=T_{11}|_{E(B_2)H}$ is the zero operator, so that the operator
decomposition
$$T_{11}=T_{11}^0\oplus T_{11}^1=T_{11}^0\oplus 0$$
is meaningless.

\medskip

\textbf{Queries 1.} In Case \ref{23c2}, the author of \cite{sm22} claimed that "if the spectrum of $T_{11}$
and $T_{nn}$ consists of only one point, then $T_{11}$ and $T_{nn}$ are a multiple of the identity."
The author of \cite{sm22} did not present the reason to show that this claim is true or indicate the
corresponding references (for example, see a certain page of a certain book [*] or see a certain paper [*]).

Indeed, it is easy to see that the spectrum of any multiple $\alpha I$ of the identity $I$ on a
Banach space $X$ consists of only one point $\{\alpha\}$. In fact, if $T=\alpha I$, it is clear that
$\alpha$ is an eigenvalue of $T$, so $\alpha\in \sigma(T)$. Moreover, when $\lambda\neq\alpha$,
the operator $T-\lambda I=(\alpha-\lambda)I$ is a bijection from $X$ onto $X$. Therefore
$\lambda\in \rho(T)=\mathbb{C}\setminus \sigma(T)$ whenever $\lambda\neq\alpha$. we conclude
$\sigma(T)=\sigma(\alpha I)=\{\alpha\}$. Perhaps it is for this reason that the author
of \cite{sm22} presented the above claim. But it seems that the above claim is not true.
Indeed, we can give a counterexample as follows in Hilbert space $l_2$, in which there is an
operator $T_{11}$ whose spectrum consists of only one point, but $T_{11}$ is not any multiple of the
identity $I$ on $l_2$.

\medskip

\textbf{Counterexample 2.} Let $T_{11}$ be an operator on $l_2$ whose representation matrix is as follows:
\begin{eqnarray*}
\left(
\begin{array}{cccc}
0&0&0&\cdots\\
1&0&0&\cdots\\
0&0&0&\cdots\\
\vdots&\vdots&\vdots&\ddots
\end{array}
\right).
\end{eqnarray*}
Then for any $x=(x_1, x_2,\cdots)\in l_2$, we have
\begin{eqnarray*}
T_{11}x=\left(
\begin{array}{cccc}
0&0&0&\cdots\\
1&0&0&\cdots\\
0&0&0&\cdots\\
\vdots&\vdots&\vdots&\ddots
\end{array}
\right)
\left(\begin{array}{c}
x_1\\
x_2\\
x_3\\
\vdots
\end{array}
\right)
=\left(\begin{array}{c}
0\\
x_1\\
0\\
\vdots
\end{array}
\right).
\end{eqnarray*}
Therefore, it is easy to see that $T_{11}$ is a  bounded linear operator on $l_2$. Moreover, as in
Part (b), it can be shown that $\sigma(T_{11})=\{0\}$. But it is clear that $T_{11}$ is not any multiple of
the identity $I$ on $l_2$.

By the way, if $T_{11}$ is a normal operator and its spectrum $\sigma(T_{11})$ consists of only one point,
then by the spectral theorem, $A_{11}$ is a multiple of the identity (cf. \cite{hrp82} p.211).
However, for M-hyponormal operator, it seems that there is no such result in the existing literatures,
and we give the above counterexample.

\end{document}